\documentstyle{amsppt}
\magnification=1200
%\UseAMSsymbols
\NoBlackBoxes
\topmatter
\title Subspaces containing biorthogonal functionals of bases of
different types
\endtitle
\rightheadtext{Biorthogonal functionals}
\author M.I.Ostrovskii
\endauthor
\address Mathematical Division, Institute for Low Temperature Physics
and Engineering, 47 Lenin avenue, 310164 Kharkov, UKRAINE \endaddress
\email mostrovskii\@ilt.kharkov.ua \endemail
\abstract
The paper is devoted to two particular cases of the following general
problem. Let $\alpha$ and $\beta$ be two types of bases in Banach
spaces.
Let a Banach space $X$ has bases of both types and a subspace
$M\subset X^*$ contains the sequence of biorthogonal functionals
of some $\alpha$-basis in $X$. Does $M$ contain a sequence of
biorthogonal
functionals of some $\beta$-basis in $X$?

The following particular cases are considered:

$(\alpha, \beta)$=(Schauder bases, unconditional bases),

$(\alpha, \beta)$=(Nonlinear operational bases, linear operational
bases).

The paper contains an investigation of some of the spaces constructed
by
S.Belle\-not in ``The $J$-sum of Banach spaces'', J. Funct. Anal. {\bf
48}
(1982), 95--106. (These spaces are used in some examples.)
\endabstract
\thanks This is an English translation of the paper published in
``Teoriya
Funktsii, Funktsional'nyi Analiz i ikh Prilozheniya'', vol. 57 (1992),
pp.~115--127
\endthanks
\endtopmatter

\document

We use the standard Banach space notation as can  be  found  in
[LT2], [PP], [S2].

1. Definition 1. Let $X$ be a Banach space with (unconditional) basis.
A  subspace
$M\subset X^*$ is called ({\it unconditionally}) {\it basic}
if it contains all biorthogonal functionals  of
some (unconditional) basis of $X$.

Basic subspaces have been studied in [DK], [O2].

{\bf Theorem 1.}  {\it Let $X$  be  a  non-reflexive  Banach  space
with  an
unconditional basis. There exists a subspace of $X^{*}$  which  is
basic
but is not unconditionally basic.}

Proof. Let $(x_{i})^{\infty }_{i=1}$ be an unconditional basis of
$X$ and $x^{*}_{i} (i\in {\Bbb N})$
be its biorthogonal functionals. Then either $(x_{i})$ is boundedly
complete or it is not. Suppose first that $(x_{i})$ is boundedly
complete.
Then $X$ can be equivalently renormed to become the dual of the space
$N=[x^{*}_{i}]^{\infty }_{i=1}$, in natural duality. (We use square
brackets
to denote the closure of linear span.)
The space $N$ is a non-reflexive Banach space and $(x^{*}_{i})$ is an
unconditional shrinking basis of it. So by James' theorems [LT2, p.~9,
22] this basis is not boundedly complete and the space $N$ contains a
sequence of blocks
$$
m_{i}=\sum^{n_i}_{k=n_{i-1}+1}a_{k}x^{*}_{k};\ || m_{i}|| =1,
$$
equivalent to the unit vector basis of $c_{0}$. Let
$$
m^{*}_{i}=\sum^{n_i}_{n_{i-1}+1}b_{k}x_{k}
$$
be chosen so that $m^{*}_{i}(m_{i})=1$ and
$\sup _{i}|| m^{*}_{i}|| =C<\infty $.
Let $P:N\rightarrow N$ be defined by
$$
P(m)=\sum^{\infty }_{i=1}m^{*}_{i}(m)m_{i}.
$$
It is clear that $P$ is a projection onto the isomorphic copy of
$c_{0}$. One can modify
Zippin's arguments [Z, p. 76] to construct shrinking basis
$(q_{i})^{\infty }_{i=1}$ in $\ker P$. We have
$$
X=([q_{i}]^{\infty }_{i=1})^{*}\oplus ([m_{i}]^{\infty }_{i=1})^{*}.
$$

It is clear that the second space is an isomorphic copy of $l_{1}$.

By $c(\alpha )$ we denote the space of all continuous  functions
on the set of all ordinals not greater than $\alpha $  provided  with
order
topology. For countable ordinal $\alpha $
we have $(c(\alpha ))^{*}=l_{1}$ and $c(\alpha )$ has  a
shrinking basis ([LT1, p. 177, 213], [LT2, p. 10]).

Let $\{s_{i}\}_{i=1}^\infty$ be a shrinking basis
of the space $L:=c(\omega ^{\omega^2} )$ and let
$(s^{*}_{i})$ be the sequence of its biorthogonal functionals. The
system
$$
\{q^{*}_{i}\}_{i=1}^\infty\cup \{s^{*}_{i}\}_{i=1}^\infty
\eqno{(1)}$$
after any enumeration preserving the order in each of the  sequences
forms a boundedly complete basis of $X$ [LT2, p. 9].
Let $M\subset X^{*}$  be  the
closure of the linear span of the biorthogonal  functionals  of  the
system (1). Since the basis (1) is boundedly  complete,  it  follows
that $M$ does not contain any proper  closed  total  subspace.  It  is
clear that $M$ is isomorphic to $L\oplus [q_{i}]$.

It is clear that $M$ is a basic subspace. We shall prove that $M$ is
not
unconditionally basic.
Let us suppose that it is not the case and let $\{u_{i}\}_{i=1}^\infty$
be an unconditional basis of $X$ whose biorthogonal functionals
$\{u^{*}_{i}\}_{i=1}^\infty$ belong to $M$.
>From the remark above we obtain that $M=[u^{*}_{i}]^{\infty }_{i=1}$
and hence $M$  have
an unconditional basis. But by Maurey - Rosenthal theorem [MR], $L$ and
therefore $M$  contains  weakly  null  normalized  sequences  with  no
unconditional subsequence, a contradiction [LT2, p.~7, 19]. Thus, if
$\{x_{i}\}_{i=1}^\infty$ is boundedly complete, then we are done.

Suppose now that $\{x_{i}\}_{i=1}^\infty$ is not boundedly complete. As
before let
us introduce $N=[x^{*}_{i}]^{\infty }_{i=1}$.
It is  easy  to  see  that  there  exists  a
functional $x^{**}\in X^{**}\backslash X$ such that $x^{**}\mid
_{N}\neq 0.$
The space $M=\ker  x^{**}\cap N$ is a
total subspace of $X^{*}$. By [DK, Theorem 3] $M$ is a basic subspace.
By [O3, Theorem 1] $M$ is not unconditionally basic.

2. Definition 2. A subspace $M\subset X^{*}$ is said to be {\it
norming}
if  there  exists
$c>0$ such that
$$
(\forall x\in X)(\sup _{0\neq f\in M}|f(x)|/||f||\ge c||x||).
$$

{\bf Remark.} M.I.Kadets [K] proved that if $X$ is separable and
$M\subset X^{*}$ is
a norming subspace, then $X$ has a nonlinear operational basis  all
of whose biorthogonal functionals are in $M$.  V.P.Fonf  [F]  proved
that every subspace with the last property is norming.

Definition 3. A subspace $M\subset X^{*}$ is said  to  be  {\it
quasibasic}  if
there exists a sequence of  continuous  finite-dimensional  linear
operators $v_{n}:X\rightarrow X (n\in {\Bbb N})$ such that
$$
1)(\forall x\in X)(\lim_{n\rightarrow \infty }||v_{n}(x)-x||=0);
$$
$$
2)(\forall n\in {\Bbb N})(\hbox{im}(v^{*}_{n})\subset M),
$$
where the operators $v^{*}_{n}$ are adjoint to $v_{n}$.

{\bf Remark.} It is easy to see that a subspace $M\subset X^{*}$ is
quasibasic  if
and only if $M$ contains all biorthogonal functionals of some linear
operational basis of $X$.
\par

Definition 4. A Banach  space $X$  is  said  to  have  the  {\it total
property of bounded approximation} (TPBA in short) if every norming
subspace $M$ of $X^{*}$ is quasibasic.
\par
This property was introduced independently  and  almost  in  the
same   time   by    I.Singer    [S1],    F.S.Vakher    [V1]    and
V.A.Vinokurov-A.N.Plichko [ViP] (we would like to note that [S1] is
based on the lecture given in 1975). Later on  this  property  was
investigated by many authors (see [G], [GP], [MP], [O1],
[V2], [VP], [VGP]),  some  of
these results are discussed in [S2, pp. 776-779,  865].  The  term
TPBA appeared in [V2]. The purpose of the present paper is to make
some additions to abovementioned works.
\par
It is clear that if $X\in $TPBA then $X$  is  separable  and  has  the
bounded  approximation  property  (BAP).  Our  aim  is   to   find
conditions under which the converse is also true.
\par
Definition 5. Let $X(1)$ and $X(2)$ be finite-dimensional  subspaces
of a Banach space $X$, such that $X(1)\subset X(2)\subset X$ and let
 $\lambda >0.$  The  pair
$(X(1),X(2))$ is  said  to  be $\lambda ${\it -approximable}  if
there  exists  a
continuous linear  operator $u:X\rightarrow X(2)$  satisfying  the
conditions
$||u||\le \lambda $ and $u|_{X(1)}=I_{X(1)}$. A sequence
$$(X(1,i),X(2,i))^{\infty }_{i=1}$$
of pairs of subspaces of $X$ is said to be {\it uniformly
approximable}
if  there
exists $0<\lambda <\infty $ such  that  all  of  the  pairs
$(X(1,i),X(2,i))$  are
$\lambda $-approximable.
\par
Definition 6. Let $U$ and $V$ be subspaces of a Banach space  $X$. The
number
$$
\delta (U,V)=\inf \{||u-v||: u\in S(U), v\in V\}
$$
is called the {\it inclination}  of $U$ to $V$.

Let $M$ be a subspace of $X^{*}$.  We  shall  denote  by $M^{\bot }$
the  set
$\{x^{**}\in X^{**}:\
(\forall x^{*}\in M)(x^{**}(x^{*})=0)\}$. It is known [PP, p. 32] that
$M$  is
norming if  and  only  if $\delta (M^{\bot },X)>0$. (We  identify $X$
with  its
canonical image in $X^{**})$.
\par
Let $\phi :X^{**}\rightarrow X^{**}/M^{\bot }$ be the  natural
quotient  mapping.
The  space
$X^{**}/M^{\bot}$ is naturally isometric to $M^{*}$. If $M$ is  a
norming  subspace
then $\phi|_{X}$ is an isomorphic embedding.

{\bf Theorem 2.} {\it Let $X$ be a separable Banach space (SBS) with
the BAP.
Let $M$ be a norming subspace of $X^{*}$. Subspace $M$ is quasibasic if
and
only  if  the  sequence $(\phi (X(1,i)),\phi (X(2,i)))^{\infty
}_{i=1}$   is   uniformly
approximable in $M^{*}$ for every uniformly approximable in $X$
sequence}
$(X(1,i),X(2,i))^{\infty }_{i=1}$.

Proof. Necessity. Let $M$ be a quasibasic subspace of $X^{*}$  and  let
$$(X(1,i),X(2,i))^{\infty }_{i=1}$$
 be a uniformly approximable sequence in $X$.  Let
$\{v_{n}\}$ be a sequence  of  operators  for  which  the  conditions
of
Definition 3 are satisfied. By Banach-Steinhaus  theorem  we  have
$\sup _{n}||v_{n}||=\beta <\infty $. Therefore we can select a
subsequence $\{v_{n(i)}\}^{\infty }_{i=1}$  of
$\{v_{n}\}$ such that
$$
(\forall x\in X(1,i))(||v_{n(i)}(x)-x||<||v_{n(i)}(x)||/\dim
(X(1,i))).
$$
Using standard reasoning (see [JRZ]) we can find operators
$A_{i}:X\rightarrow X
\ (i\in {\Bbb N})$ such  that $||A_{i}||\le 2$ and
$$
(\forall x\in X(1,i))(A_{i}v_{n(i)}(x)=x).
$$
Since the sequence $(X(1,i),X(2,i))^{\infty }_{i=1}$ is uniformly
approximable,
then for some $\lambda <\infty $
there exists a sequence $\{u_{i}\}^{\infty }_{i=1}$  of  operators,
$u_{i}:X\rightarrow X$ such that
$$
(\forall i\in {\Bbb N})(\hbox{im}(u_{i})\subset X(2,i));
$$
$$
\sup _{i}||u_{i}||=\lambda <\infty ;
$$
$$
u_{i}|_{X(1,i)}=I_{X(1,i)}.
$$
Let $T_{i}=u_{i}A_{i}v_{n(i)}$. We have
$$
T_{i}|_{X(1,i)}=I_{X(1,i)};
\eqno{(2)}$$
$$
\hbox{im}(T_{i})\subset X(2,i);
\eqno{(3)}$$
$$
\hbox{im}(T^{*}_{i})=\hbox{im}(v^{*}_{n(i)}A^{*}_{i}u^{*}_{i})\subset
\hbox{im}(v^{*}_{n(i)})\subset M;
\eqno{(4)}$$
$$
||T_{i}||\le 2\lambda \beta .
\eqno{(5)}$$
Conditions (3) and (4) means that $T_{i}$ can be represented  in  the
form $T_{i}(x)=\sum^{n(i)}_{k=1}f^{i}_{k}(x)x^{i}_{k}$,
where $x^{i}_{k}\in X(2,i),\ f^{i}_{k}\in M$. Let operators
$R_{i}:M^{*}\rightarrow M^{*}$ be given by the equalities
$R_{i}(m^{*})=\sum^{n(i)}_{k=1}m^{*}(f^{i}_{k})\phi (x^{i}_{k})\ (i\in
{\Bbb N})$.
It is clear that $R_{i}$ are $\sigma (M^{*},M)$-continuous and that
$\phi (B(X))$ (where
$B(X)$ is the closed unit ball of $X$) is $\sigma (M^{*},M)$-dense in
some  ball
of non-zero radius of $M^{*}$ [PP, p.~32].  By  (5)  it  follows  that
$R_{i}$ are uniformly continuous operators on $M^{*}$. By (2) it
follows that
$R_{i}|_{\phi X(1,i)}=I_{\phi X(1,i)}$, and by (3) it  follows  that
$\hbox{im}(R_{i})\subset \phi X(2,i)$.
The necessity is proved.

Sufficiency. If $X$ has the BAP and is separable then it  is  easy
to  find  a  sequence $\{X(i)\}^{\infty }_{i=1}$  of  subspaces  of
$X$  such  that
$X(1)\subset X(2)\subset \ldots
\subset X(n)\subset \ldots$;
cl$(\cup ^{\infty }_{n=1}X(n))=X$    and     the     pairs
$(X(1,i),X(2,i))=(X(i),X(i+1))$  forms  a   uniformly   approximable
sequence.   Our   supposition   implies    that    the    sequence
$(\phi (X(1,i)),\phi (X(2,i)))$ is a uniformly approximable sequence in
$M^{*}$.
Let the operators $R_{i}:M^{*}\rightarrow M^{*}$ be such that $\sup
_{i}||R_{i}||<\infty ; \hbox{im}(R_{i})\subset \phi X(2,i)$
and $R_{i}|_{\phi X(1,i)}=I_{\phi X(1,i)}$.

{\bf Lemma 1} [JRZ, p. 494]. {\it Let $L$ and $N$ be Banach spaces
with $\dim (N)<\infty $.
Let $F$ be a finite dimensional subspace of $L^{*}$, let $Q$ be an
operator
from $L^{*}$ into $N$ and let $\varepsilon >0.$
Then  there  is  a  weak$^{*}$  continuous
operator $R$ from $L^{*}$ to $N$ such that $R|_{F}=Q|_{F}$ and}
$||R||\le ||Q||(1+\varepsilon )$.

By this lemma we may without  loss  of  generality  assume  that
operators $R_{i}$ are weak$^{*}$ continuous, i.e.
$$
R_{i}(m^{*})=\sum^{n(i)}_{k=1}m^{*}(f^{i}_{k})\phi x^{i}_{k},
$$
where $f^{i}_{k}\in M,\ x^{i}_{k}\in X(2,i)$.
Let  operators $T_{i}:X\rightarrow X$   be   given   by
$T_{i}(x)=\sum^{n(i)}_{k=1}f^{i}_{k}(x)x^{i}_{k}$. We have
$$
\sup _{i}||T_{i}||<\infty ;
\eqno{(6)}$$
$$
T_{i}|_{X(1,i)}=I_{X(1,i)};
\eqno{(7)}$$
$$
\hbox{im}(T^{*}_{i})\subset M
\eqno{(8)}$$
By (6), (7) and the equality cl$(\cup ^{\infty }_{n=1}X(n))=X$ we
obtain:
$$
(\forall x\in X)(\lim_{n\rightarrow \infty }||T_{n}(x)-x||=0).
$$
By (8) it follows that $M$ is quasibasic. The theorem is proved.
\bigskip

Using this theorem we can obtain the following result of [MP].

{\bf Corollary.} {\it Let $X$ be a SBS with the  BAP.  Let $M$  be  a
norming
subspace of $X^{*}$, such that the
subspace $M^{\bot }\subset X^{**}$  has  a  complement,
which  contains $X$. Then $M$ is quasibasic.}

Proof. Let us show that $M$ satisfies the assumptions  of  theorem
2. Let
$$(X(1,i),X(2,i))^{\infty }_{i=1}$$
be  a uniformly approximable sequence in
$X$. Let $Y$ be a complement of $M^{\bot }$,
such that $Y\supset X$. It  is  clear  that
the restriction of $\phi $ to $Y$ is an  isomorphism  between $Y$  and
$M^{*}$.
Therefore,  it  is  sufficient   to   show   that   the   sequence
$(X(1,i),X(2,i))^{\infty }_{i=1}$ is uniformly approximable in $Y$.
But it is clear
that  the  second  conjugates   of   operators,   which   unifomly
approximate  pairs $(X(1,i),X(2,i))$   in $X$,   realize   uniform
approximation of pairs $(X(1,i),X(2,i))$ in $X^{**}$ and, hence, in
$Y$.

{\bf Remarks.}  1.  Existence  of  the  complement  mentioned  in  the
corollary is not necessary. It follows from the  following  result
of [VP]: every ${\Cal L}_{\infty }$-space (in the sense of
Lindenstrauss-Pe\l czy\'nski)
has the TPBA.

2. The assertion of the corollary becomes wrong if we  omit  the
condition $Y\supset X$ (see Remark after Theorem 3).

Theorem 2 reduces the problem of characterization of the TPBA to
the following one: for what SBS with the BAP there exist  a  weak$^{*}$
closed subspace $H$ of $X^{**}$ such  that $\delta (H,X)>0,$  and  the
quotient
mapping $Q:X^{**}\rightarrow X^{**}/H$ maps
some uniformly approximable sequence in $X$
on  the sequence which is not uniformly approximable in $X^{**}/H$.  We
shall describe one of the approaches to this problem.

Definition 7. Let $f:{\Bbb N}\rightarrow (0,+\infty )$. We shall say
that
a sequence
$$(Z(1,i),Z(2,i))^{\infty }_{i=1}\ (Z(1,i)\subset Z(2,i))$$
of pairs
of subspaces  of  a  Banach
space $Z$ is $f$-{\it approximable} if there  exists  a
sequence $\{u_{i}\}^{\infty }_{i=1}$ of
linear    continuous    operators $u_{i}:Z\rightarrow Z(2,i)$
such     that
$u_{i}|_{Z(1,i)}=I_{Z(1,i)}$ and $\sup _{i}(||u_{i}||/f(i))<\infty $.

{\bf Proposition.} {\it Let $H$ be a weak$^{*}$ closed subspace
of $X^{**}$  and  let
$\delta (H,X)>0.$ Let $\phi $
denote the quotient mapping $\phi :X^{**}\rightarrow X^{**}/H$.  Let
us
suppose  that $X$  contains  a  uniformly   approximable   sequence
$(X(1,i),X(2,i))^{\infty }_{i=1}$ such that for some sequence
$(Y(1,i),Y(2,i))^{\infty }_{i=1}$
of  pairs  of  subspaces  of $X^{**}$  the  following  conditions  are
satisfied.
$$
\phi X(1,i)=\phi Y(1,i);\ \phi X(2,i)=\phi Y(2,i);
\eqno{(9)}$$
$$
(\forall i\in {\Bbb N})(\delta (Y(2,i),H)>0);
\eqno{(10)}$$
and,  furthermore,  the   sequence $(Y(1,i),Y(2,i))^{\infty }_{i=1}$
is   not
$f$-approximable in $X^{**}$ for $f$ defined  by $f(i)=1/\delta
(Y(2,i),H)$.  Then
subspace $H^\top$ (where $H^\top=\{x^{*}\in X^{*}:\
(\forall x^{**}\in H)(x^{**}(x^{*})=0)\})$  is  a  norming
nonquasibasic subspace.}

Proof. Let us suppose that it is not the case and apply  Theorem
2. We obtain that the sequence $(\phi X(1,i), \phi X(2,i))$  is
uniformly
approximable in $X^{**}/H$. This means
that for some $0<\lambda <\infty $ there  exist
operators $u_{i}:X^{**}/H\rightarrow \phi X(2,i)$ such that
$$
u_{i}|_{\phi X(1,i)}=I_{\phi X(1,i)}
\eqno{(11)}$$
and $||u_{i}||\le \lambda $.
Let  operators $v_{i}:X^{**}\rightarrow Y(2,i)$  be  defined  by
$$v_{i}=(\phi|_{Y(2,i)})^{-1}u_{i}\phi .$$

This  operators  are   well-defined   because
$\hbox{im}(u_{i})\subset \phi X(2,i)=\phi Y(2,i)$, and
the inequality $\delta (Y(2,i),H)>0$  implies
that the inverse of $\phi|_{Y(2,i)}$  exists.  It  is  easy  to  see
that
$||(\phi|_{Y(2,i)})^{-1}||=f(i)$. Therefore,
$$
(\forall i\in {\Bbb N})(||v_{i}||\le \lambda f(i)).
$$
Furthermore, by (9) and (11)  we  have $v_{i}|_{Y(1,i)}=I_{Y(1,i)}$.
This
contradicts  the  assumption  that $(Y(1,i),Y(2,i))^{\infty }_{i=1}$
is   not
$f$-approximable. The proposition is proved.
\par
The verification  of  the  conditions  of  the  proposition  for
concrete spaces is laborious. Therefore, the  following  criterion
is of interest.

{\bf Theorem 3.} {\it Let $X^{**}$ contains a reflexive uncomplemented
subspace
$Y$ which is isomorphic to a complemented subspace $Z$  of $X$  and  is
such that $\delta (Y,X)>0.$ Then $X$ does not have the TPBA.}

Proof. Let $T:Y\rightarrow Z$ be an isomorphism. Let us consider the
subspace
$H=\{y-Ty: y\in Y\}$ of $X^{**}$. We shall check that  it  satisfies
all  the
conditions of the proposition with $f(i)\equiv C>0.$

Since $\delta (X,Y)>0,$ then $H$ is isomorphic to $Y$ and, hence,
reflexive.
Therefore, subspace $H$ is weak$^{*}$ closed by Krein-Smulian theorem.
It
is easy to see  that $\delta (H,X)>0$  and,  therefore,  [PP,
p.~29--34]
subspace $M=H^\top\subset X^{*}$ is norming.

It is clear that we may restrict ourselves to the case when $X$ is
a SBS with the BAP. In this case $Z$ is also a SBS with the BAP. Let
$Z(1)\subset Z(2)\subset \ldots
\subset Z(n)\subset \ldots
$  be  a   sequence   of   finite-dimensional
subspaces of $Z$ such that
$$
\hbox{cl}(\cup ^{\infty }_{n=1}Z(n))=Z,
\eqno{(12)}$$
and the sequence $(Z(i),Z(i+1))$ is uniformly approximable in $Z$ and,
hence, is uniformly approximable in $X$.

Let  us  introduce  the  following   sequences   of   subspaces:
$X(1,i)=Z(i),\ X(2,i)=Z(i+1),\ Y(1,i)=T^{-1}Z(i),\ Y(2,i)=T^{-1}Z(i+1)$.
\par
Let  us  show  that  the  sequence $(Y(1,i),Y(2,i))^{\infty }_{i=1}$
is   not
uniformly approximable in $X^{**}$. In fact, if we assume that for some
operators $u_{i}:X^{**}\rightarrow Y(2,i)$ we have $\sup
_{i}||u_{i}||<\infty $ and
$$
u_{i}|_{Y(1,i)}=I_{Y(1,i)},
\eqno{(13)}$$
then by reflexivity of $Y$ we can define the operator
$u:X^{**}\rightarrow Y$ by the
equality $u(x)=w-\lim_{A}u_{i}(x)$, where $A$ is some ultrafilter  on
${\Bbb N}$.  By
(12)  and  (13)  this  operator  is  a  projection  onto  $Y$.  This
contradicts the fact that $Y$ is uncomplemented. It is easy to check
that  all  the  other  conditions  of  the  proposition  are  also
satisfied. The theorem is proved.

{\bf Corollary.} {\it There exists a SBS $X$ with a basis which is
isometric
to its bidual but does not have the TPBA.}

Proof. Let $X=(\sum^{\infty }_{i=1}\oplus J)_{p}\ (p\neq 1,2,\infty )$,
where $J$ is  the  James'  space
(nonreflexive space with a basis, such that $J^{**}$ is isometric to
$J$,
and $J$ has codimension one in $J^{**}$ (see [LT2, p.~25])). It  is
clear
that $X$  has  a  basis  and  is  isometric  to  its  second  dual.
Furthermore, we have $X^{**}=X\oplus l_{p}$. By well-known results
([BDGJN], [R])
$l_{p}$ contains an uncomplemented subspace isomorphic to $l_{p}$.  On
the
other hand, $X$ contains a complemented subspace isomorphic  to
$l_{p}$.
We are in the conditions of Theorem 3.

{\bf Remark.} If we develop the construction  of  Theorem  3  for  the
space $X$  from  the  corollary,  then  the  subspace $H$  would  be
complemented in $X^{**}$.

In fact,  let $P:X\rightarrow Z$  be  the  projection,  whose
existence  is
supposed and let $Q:X^{**}\rightarrow X$ be the projection
corresponding  to  the
decomposition
$X^{**}=X\oplus l_{p}$. Then $PQ$ is a projection of $X^{**}$  on $Z$
and
$PQ|_{Y}=0.$ Therefore, the operator $(I_{X^{**}}-T^{-1})PQ$ is  a
projection  of
$X^{**}$ onto $H$.

It turns out that a SBS $X$ with the BAP but without the TPBA need
not satisfy the conditions of Theorem 3.

{\bf Theorem 4.} {\it There exists a SBS $X$ with a basis
such  that $X^{**}=X\oplus Y$
and $Y$ does not contain infinite-dimensional  subspaces  which  are
isomorphic to subspaces of $X$, but $X\not\in $TPBA.}

Proof.  We  need  to  use  the   variant   of   the   proof   of
James-Lindenstrauss  theorem  ([J],  [L]),  which   is   due   to
S.F.Bellenot [B].  We use the following  particular  case  of  the
construction of [B].

Let $(X_{n})^{\infty }_{n=0}$ be  an  increasing  sequence  of
finite-dimensional
subspaces of a Banach space $Z$, such that
cl$(\cup ^{\infty }_{n=0}X_{n})=Z$. For ease of
notation we adopt the convention that $X_{0}=\{0\}$.
\par
Let $(x_{i})^{\infty }_{i=0}$
be a  sequence  with $x_{i}\in X_{i}$.  If $(x_{i})$  is  finitely
non-zero, then we define the norm $||\cdot ||_{J}$ by
$$
2||(x_{i})^{\infty }_{i=0}||^{2}_{J}=
\sup (\sum^{k-1}_{i=1}||x_{p(i)}-x_{p(i+1)}||^{2}+||x_{p(k)}||^{2}),
$$
where the $\sup $  is  over  all  integer  sequences
$(p(i))^{k}_{i=1}$  with
$0\le p(1)<p(2)<\ldots
<p(k)$. The completion of the space of all  finitely
non-zero sequences under this norm will be denoted by $J(X_{n})$.
\par
We shall call the sequence
$(x_{i})^{\infty }_{i=0}, x_{i}\in X_{i}$  {\it eventually  constant},
if for some $n\in {\Bbb N}$ we have $x_{n}=x_{n+1}=x_{n+2}=\ldots
$. We endow  the  space  of
all eventually constant sequences with the semi-norm
$$
||(x_{i})^{\infty }_{i=0}||_{\Omega }=\lim_{k\rightarrow \infty
}||x_{k}||,
$$
and denote this space by $\Omega (X_{n})$. We denote by $K(X_{n})$
the  space  of
all sequences $(x_{i})$ with $x_{i}\in X_{i}$ and whose norm
$$
||(x_{i})^{\infty }_{i=0}||_{K}=\sup _{n}||(x_{0},\ldots
,x_{n},0,\ldots
,0,\ldots
)||_{J}
$$
is finite. It is clear that $K(X_{n})\supset \Omega (X_{n})$.

{\bf Theorem 5} [B]. {\it Let the sequence
$(X_{n})^{\infty }_{n=0}$ be as above.
Then\par
(I) $\Omega (X_{n})$ is dense in $K(X_{n})$.

(II) $(J(X_{n}))^{**}=K(X_{n})$ and $(J(X_{n}))^{**}/J(X_{n})$ is
isometric to $Z$.

(III) If the  spaces $X_{n}\ (n\in {\Bbb N})$  have  uniformly
bounded  basic
constants, then $J(X_{n})$ has a basis.}

Let us turn to the proof of theorem 4. We fix some $1<p<2$ and let
$Z=l_{p}$. For $X_{n}$ we take the linear spans of the first $n$
vectors  of
the  unit  vector  basis  of $l_{p}$.  Let  us  introduce  the   space
$X=J(X_{n})\oplus l_{2}$. The space $X$ has a basis by part III of
Theorem 5.

{\bf Lemma 2.} {\it The space $X^{**}$ can be represented
in the form: $X^{**}=X\oplus l_{p}$.}

Proof. Denote by $\{e^{n}_{i}\}^{n}_{i=1}$ the unit
vector basis  of $X_{n}$.  Let  us
introduce the vectors
$$
f_{i}=(0,\ldots
,0,e^{i}_{i},e^{i+1}_{i},\ldots
)\in K(X_{n}).
$$
Let us show that the sequence $\{f_{i}\}^{\infty }_{i=1}$ is equivalent
to the  unit
vector basis of $l_{p}$. We have
$$
||\sum^{\infty
}_{i=1}a_{i}f_{i}||_{K}=||(0,a_{1}e^{1}_{1},a_{1}e^{2}_{1}+a_{2}e^{2}_{2},\ldots
,\sum^{n}_{i=1}a_{i}e^{n}_{i},\ldots
)||_{K}.
$$
Recall the definition of $K$-norm and choose $p(1)=0$ and $p(2)=n$.  We
obtain
$$
||\sum^{\infty }_{i=1}a_{i}f_{i}||_{K}\ge
||\sum^{n}_{i=1}a_{i}e^{n}_{i}||=(\sum^{n}_{i=1}|a_{i}|^{p})^{1/p}.
$$
Since this inequality is valid for every $n\in {\Bbb N}$, then we have
$$
||\sum^{\infty }_{i=1}a_{i}f_{i}||_{K}\ge
(\sum^{\infty}_{i=1}|a_{i}|^{p})^{1/p}.
$$
On the other hand, we have
$$
||\sum a_{i}f_{i}||=
2^{-1/2}\sup_{(p(i))}(\sum^{k-1}_{i=1}
(\sum^{p(i+1)}_{s=p(i)+1}|a_{s}|^{p})^{2/p}+
$$
$$
(\sum^{p(k)}_{s=1}|a_{s}|^{p})^{2/p})^{1/2}\le
2^{-1/2}\sup
_{(p(i))}(\sum^{k-1}_{i=1}\sum^{p(i+1)}_{s=p(i)+1}|a_{s}|^{p}+
$$
$$
\sum^{p(k)}_{s=1}|a_{s}|^{p})^{1/p}<
2^{1/2}\sup _{(p(i))}(\sum^{p(k)}_{s=1}|a_{s}|^{p})^{1/p}=
2^{1/2}(\sum^{\infty }_{s=1}|a_{s}|^{p})^{1/p}.
$$
>From here and from the proof of Theorem 5 in [B] it follows that
the restriction of the quotient mapping $K(X_{n})\rightarrow
K(X_{n})/J(X_{n})$
to  the
closure of the linear span of $\{f_{i}\}$ is an isomorphism. The lemma
is
proved.
\par
{\bf Lemma 3.} {\it Every infinite-dimensional subspace  of $X$
contains  a
subspace isomorphic to} $l_{2}$.

Proof. Since $X=J(X_{n})\oplus l_{2}$, then it  is  sufficient  to
show  that
every infinite-dimensional subspace of $J(X_{n})$ contains  a  subspace
isomorphic to $l_{2}$.
\par
It is easy to see (it is shown in  the  proof  of  part  III  of
Theorem 5 in [B]) that the vectors $f^{n}_{i}=(0,\ldots
,0,e^{n}_{i},0,\ldots
)$  form  a
basis of $J(X_{n})$.
\par
The  equality $(J(X_{n}))^{**}=J(X_{n})\oplus l_{p}$   implies
separability   of
$(J(X_{n}))^{*}$. Therefore, every infinite-dimensional subspace of
$J(X_{n})$
contains a weakly null sequence $(x_{k})^{\infty }_{k=1}$ which is
bounded away from
zero. By the well-known arguments [LT2, p.~7] it follows that we can
select a subsequence $(x_{n(k)})^{\infty }_{k=1}$
of $(x_{k})$, which  is  equivalent  to
the sequence of the form
$$
h_{k}=\sum^{r(k+1)-1}_{n=r(k)+1}(\sum^{n}_{i=1}a^{n}_{i}f^{n}_{i}).
$$
It can be directly verified that the sequence $(h_{k})$ is equivalent
to the unit vector basis of $l_{2}$. The lemma is proved.

{\bf Lemma 4.} $X\not\in ${\it TPBA.}

Proof. It is known [BDGJN] that  for  every $1<p<2$  there  exists  a
sequence $\{W_{i}\}^{\infty }_{i=1}$ of finite-dimensional subspaces
of $l_{p}$  such  that
the following conditions are satisfied: $\dim (W_{i})=i$;
$$
\sup _{i}d(W_{i},l^{i}_{2})=C<\infty ;
$$
$$
(\exists 0<c_{1}\le c_{2}<\infty )(\forall \{w_{i}\}^{\infty
}_{i=1};\ w_{i}\in W_{i})
$$
$$
(c_{1}(\sum ||w_{i}||^{p})^{1/p}\le ||\sum w_{i}||\le c_{2}(\sum
||w_{i}||^{p})^{1/p});
$$
and the sequence $(W_{i},W_{i})^{\infty }_{i=1}$ is not uniformly
approximable.
\par
Let $W$=cl(lin$(\cup ^{\infty }_{i=1}W_{i}))$.
By Lemma 2 we have $X^{**}=X\oplus l_{p}=J(X_{n})\oplus l_{2}\oplus
l_{p}$.
Let us represent $l_{2}$ as an infinite direct sum:
$l_{2}=(\sum^{\infty }_{i=1}U_{i})_{2}$, where
$U_{i}$ are subspaces isometric to $l^{i}_{2}$.
Let  the  isomorphisms $T_{i}:W_{i}\rightarrow U_{i}\
(i\in {\Bbb N})$ are such that
$$
||T_{i}||\le 1;\ ||T^{-1}_{i}||\le C.
\eqno{(14)}$$
Let  us  introduce  the  operator $T:W\rightarrow l_{2}$   by   the
equality
$T((w_{i})^{\infty }_{i=1})=(T_{i}w_{i})^{\infty }_{i=1}$.
It is clear that $T$ is a bounded operator and
that
$$
TW_{i}=U_{i}.
\eqno{(15)}$$
Let $H=\{x-Tx: x\in W\}\subset X^{**}$.
Let us check that $H$  satisfies  all  the
conditions of the proposition.
\par
The subspace $H$ is weak$^{*}$ closed by reflexivity of $W$. It is
clear
that $\delta (H,X)>0.$ Let
$X(1,i)=X(2,i)=U_{i},\ Y(1,i)=Y(2,i)=W_{i}$. It is clear
that
$$(X(1,i),X(2,i))^{\infty }_{i=1}$$
   is    uniformly    approximable    and
$(Y(1,i),Y(2,i))^{\infty }_{i=1}$ is not. Condition (15) implies  (9),
and  (14)
implies (10). Moreover, we have $\inf _{i}\delta (Y(2,i),H)>0.$
The  lemma  is
proved.
\par
Theorem 4 follows immediately from Lemmas 2, 3 and 4.
\par
3. Definition 8. Let $X$ be a subspace of a Banach space $Z$ and let
$M$ be a subspace of $X^{*}$. The subspace $M$  is  said  to  be  {\it
boundedly
extendeable} onto $Z$ if there exists an
isomorphic embedding $\pi :M\rightarrow Z^{*}$
such that
$$
(\forall f\in M)(\forall x\in X)((\pi (f))(x)=f(x)).
$$

{\bf Theorem 6.} {\it Let $X$ be a SBS with the BAP. The  space $X$
does  not
have the TPBA if and only if there exist a Banach  space $Z$  such
that $X$ is a  subspace  of $Z$  and  the  following  conditions  are
satisfied:

(a) $X^{\bot }$ is uncomplemented in $Z^{*}$;
\par
(b) $X^{*}$  contains  a  norming  subspace $M$  which  is  boundedly
extendeable onto $Z$.}
\par
Proof. Sufficiency is proved in [V2]. Here is a shorter  proof
of it.
\par
Let us show that the subspace $M$  from  the  formulation  of  the
theorem is  not  quasibasic.  Let  us  assume  the  contrary.  Let
finite-dimensional continuous operators $u_{n}:X\rightarrow X$ be such
that
$$
(\forall x\in X)(\lim_{n\rightarrow \infty }||u_{n}(x)-x||=0);
$$
$$
(\forall n\in {\Bbb N})(\hbox{im}(u^{*}_{n})\subset M).
$$
Therefore, the operators $u_{n}$ can be represented in the  following
form: $u_{n}(x)=\sum^{p(n)}_{i=1}x^{*}_{i,n}(x)x_{i,n}$,
where $x^{*}_{i,n}\in M;\ x_{i,n}\in X$. Let us  denote  by
$r:Z^{*}\rightarrow X^{*}$  the  operator  of  the
restriction  and  by $\pi :M\rightarrow Z^{*}$  the
operator,  whose  existence  follows  from  the  definition  of  a
boundedly extendeable subspace. Let  us  introduce  the  operators
$\alpha _{n}:Z^{*}\rightarrow Z^{*}\ (n\in {\Bbb N})$ by the equalities
$\alpha _{n}(z^{*})=\sum^{p(n)}_{i=1}z^{*}(x_{i,n})
\pi (x^{*}_{i,n})=\pi u^{*}_{n}r(z^{*})\ (n\in {\Bbb N})$.
\par
It  is  easy  to  see  that  the  sequence
$\{\alpha _{n}\}^{\infty }_{n=1}$  is  uniformly
bounded.   Let   the    operator $Q:Z^{*}\rightarrow Z^{*}$    be
defined    by
$Q(z^{*})=w^{*}-\lim_{A}\alpha _{n}(z^{*})$, where $A$ is
some ultrafilter  on  the  set  of
natural numbers. Let us show that $Q$  is  a  projection  and  that
$\ker (Q)=X^{\bot }$.
\par
The relation $X^{\bot }\subset \ker (Q)$ follows immediately from  the
definition
of $Q$. Furthermore, we have
$$
rQz^{*}=w^{*}-\lim_{A}r\alpha
_{n}(z^{*})=w^{*}-\lim_{A}u_{n}r(z^{*})=r(z^{*}).
\eqno{(*)}$$
Therefore, $\ker (Q)=X^{\bot }$. The equality $Q^{2}=Q$ follows by
(*)  and  the
fact that $Q(z^{*})$  depends  only  on $r(z^{*})$.  Therefore,
$X^{\bot }$  is  a
complemented subspace of $Z^{*}$.  This  contradiction  completes  the
proof.
\par
Necessity. Let $X\not\in $TPBA and $M$ be a norming nonquasibasic
subspace
of $X^{*}$. Let $Z=M^{*}$. There is a natural isomorphic embedding of
$X$ into
$Z$. Therefore (after corresponding renorming) we may consider $X$  as
a subspace of $Z$. The subspace $M$  is  a  norming  subspace  of
$X$.
Furthermore, $M$ is boundedly extendeable to $Z$ in a natural way.  It
remains to prove that $X^{\bot}$ is uncomplemented subspace of
$Z^{*}$.  Assume
the contrary. In this case $M^{**}=Z^{*}$ can be represented in  the
form
$X^{\perp }\oplus U$, moreover $U$ is isomorphic
to $X^{*}$  in  a  natural  way.  Since
$X\in $BAP then there exists vectors
$\{x_{i,n}\}^{p(n)}_{i=1}\subset X$ and
$\{x^{*}_{i,n}\}^{p(n)}_{i=1}\subset X^{*}$  such
that
$$
(\forall x\in X)
(x=\lim_{n\rightarrow \infty }\sum^{p(n)}_{i=1}x^{*}_{i,n}(x)x_{i,n}).
$$
We denote by $S:X^{*}\rightarrow U$ the natural isomorphism.  Let  us
introduce
the sequence $\{T_{n}\}^{\infty }_{n=1}$ of the operators,
$T_{n}:Z\rightarrow Z$ by the equalities:
$$
T_{n}(z)=\sum^{p(n)}_{i=1}(Sx^{*}_{i,n})(z)x_{i,n}.
$$
This sequence is uniformly bounded. It converges to the identity
operator on $X\subset Z$. By Lemma 1 and separability of $X$ we  can
find  a
sequence of weak$^{*}$ continuous operators on $M^{*}=Z$,  such  that
their
restrictions to $X$ converge to the identity operator. Hence, $M$ is a
quasibasic subspace of $X^{*}$. The theorem is proved.

4. The result of [MP] cited after Theorem 2 implies  that  if  a
SBS $X$ with the BAP is such that every closed norming subspace $M$ of
$X^{*}$ has a finite codimension, then $X\in $TPBA. Therefore, it is
useful
to study the class of such spaces and to compare it with the class
of quasireflexive SBS with the BAP. (Recall that a Banach space $X$
is called {\it quasireflexive} if $\dim (X^{**}/X)<\infty )$.
\par
W.J.Davis and W.B.Johnson [DJ] gave examples of nonquasireflexive
SBS such that every closed norming subspace $M$ of $X^{*}$ is  of
finite
codimension. The  argument  in  [DJ]  is  based  on  the  following
observation. If $M$ is a norming subspace of $X^{*}$  then,  on  the
one
hand, $M^{\bot}\subset X^{**}$ is isomorphic to a
subspace to a  subspace  of $X^{**}/X$
and, on the other hand, $M^{\bot }$ is isometric to $(X^{*}/M)^{*}$.
Therefore,  if
$X^{**}/X$ does not contain infinite-dimensional  subspaces  which  are
isomorphic to dual spaces, then $X^{*}$ does not contain closed norming
subspaces of infinite codimension. The purpose of the  final  part
of the present paper is to show that  the  converse  statement  is
false.

{\bf Theorem 7.} {\it There exists a Banach space $Y$ with a basis,
such that
the quotient space $Y^{**}/Y$ is an infinite-dimensional reflexive SBS,
but $Y^{*}$ does not  contain  closed  norming  subspaces  of  infinite
codimension.}

Proof. We use the construction due to S.F.Bellenot [B], which is
described above. Let $p>2,\ Z=l_{p}$ and $X_{n}$ be the linear  span
of  the
first $n$ elements of the unit vector basis of $l_{p}$. Let
$Y=J(X_{n})$.

{\bf Lemma 5.} {\it The space $Y^{**}$ does not contain isomorphic
copies of $l_{p}$.}

Proof. Assume the contrary. Let $\{f_{i}\}^{\infty }_{i=1}$ be a
sequence in $Y$, that
is equivalent to the  unit  vector  basis  of $l_{p}$.  Since $Y^{**}$
is
isometric to $K(X_{n})$, then we can represent $f_{i}$ in the form

$$f_{i}=(x^{i}_{0},x^{i}_{1},\ldots
,x^{i}_{n},\ldots),
\eqno{(16)}$$
where $x^{i}_{n}\in X_{n}$. By part I  of  theorem  5  we  may  assume
that  the
sequence (16) is eventually constant. Let
$f_{1}=(x^{1}_{0},x^{1}_{1},\ldots
,x^{1}_{n(1)-1},x^{1},x^{1},\ldots
,x^{1},\ldots
)$.
\par
Since the sequence $\{f_{i}\}^{\infty }_{i=1}$ is
weakly null, then for  every $\varepsilon _{2}>0$
we can find a natural number $m(2)$ such that $f_{m(2)}$  satisfies
the
condition $||x^{m(2)}_{0}||+\ldots
+||x^{m(2)}_{n(1)+1}||<\varepsilon _{2}$. We have
$$
f_{m(2)}=(x^{m(2)}_{0},\ldots
,x^{m(2)}_{n(2)-1},x^{m(2)},\ldots
,x^{m(2)},\ldots
).
$$
It is clear that we may suppose that $n(2)-1>n(1)+1.$ Let $\varepsilon
_{3}>0.$  We
can find a natural number $m(3)$  such  that $f_{m(3)}$  safisfies  the
condition $||x^{m(3)}_{0}||+\ldots
+||x^{m(3)}_{n(2)+1}||<\varepsilon _{3}$. We have
$f_{m(3)}=(x^{m(3)}_{0},\ldots
,x^{m(3)}_{n(3)-1},x^{m(3)},\ldots
,x^{m(3)},\ldots
)$.

We continue in an obvious manner.
\par
We choose the sequence $\{\varepsilon _{k}\}^{\infty }_{k=1}$
quickly  converging  to  zero  in
order that the sequence
$$
g_{1}=f_{1};
$$
$$
g_{2}=(0,\ldots
,0,x^{m(2)}_{n(1)+2},\ldots
,x^{m(2)}_{n(2)-1},x^{m(2)},\ldots
,x^{m(2)},\ldots
);
$$
$$
g_{3}=(0,\ldots
,0,x^{m(3)}_{n(2)+2},\ldots
,x^{m(3)}_{n(3)-1},x^{m(3)},\ldots
,x^{m(3)},\ldots
);
$$
$$
\ldots
\ldots
\ldots
\ldots
\ldots
\ldots
\ldots
\ldots
\ldots
\ldots
\ldots
\ldots
\ldots
\ldots
\ldots
\ldots
\ldots
..
$$
\noindent is equivalent to the unit vector basis of $l_{p}$.  Later
on  for  the
sake of convenience we let $m(1)=1,\ n(0)=-1.$
\par
Let $\{p(k,i)\}^{\infty }_{k=1},^{s(k)}_{i=1}$
be a collection  of  natural  numbers  such
that $n(k-1)+1=p(k,1)<\ldots
<p(k,s(k))=n(k)$ and
$$
\sum^{s(k)-1}_{i=1}||x^{m(k)}_{p(k,i+1)}-x^{m(k)}_{p(k,i)}||^{2}\ge
||g_{k}||^{2}_{K}.
$$
Let us estimate $||\sum^{}_{k}a_{k}g_{k}||_{K}$. For this let us
consider the sequence
$(p(i))^{\infty }_{i=1}$ that consists of the following integers:
\par
$$
p(1,1)<p(1,2)<\ldots
<p(1,s(1))<
$$
$$
p(2,1)<p(2,2)<\ldots
<p(2,s(2))<\ldots
<
$$
$$
p(r,1)<p(r,2)<\ldots
<p(r,s(r)).
$$
We obtain
\par
$$
2||\sum^{r}_{k=1}a_{k}g_{k}||\ge
\sum^{r}_{k=1}a^{2}_{k}||g_{k}||^{2}_{K}.
$$
Since $p>2,$ then this estimate contradicts the fact that
$\{g_{k}\}^{\infty }_{k=1}$
is equivalent to the unit vector basis of $l_{p}$. The lemma is
proved.
\par
Proof of Theorem 7. The space $Y$ has a basis by the part  III  of
theorem 5. Part II of theorem 5 implies that $Y^{**}/Y$ is isometric to
$l_{p}$.
\par
Let us assume that $Y^{*}$ contains  a  closed  norming  subspace  of
infinite codimension. Then $M^{\perp }\subset Y^{**}$ is isomorphic to
a  subspace  of
$Y^{**}/Y$, i.e., to a subspace of $l_{p}$. Hence, $M^{\perp }$
contains  a  subspace
isomorphic to $l_{p}$. By Lemma 5 this is impossible.  The  theorem  is
proved.
\par
\centerline{REFERENCES
}

[B] S.F.Bellenot, The $J$-sum of Banach spaces, J.  Funct.  Anal.  48
(1982), 95-106.

[BDGJN] G.Bennett, L.E.Dor, V.Goodman, W.B.Johnson and  C.M.Newman,  On
uncomplemented subspaces  of $L_{p}, 1<p<2,$  Israel  J.  Math.  26
(1977), 178-187.

[DJ]  W.J.Davis  and  W.B.Johnson,  Basic  sequences   and   norming
subspaces in non-quasi-reflexive Banach spaces, Israel J.  Math.
14 (1973), 353-367.

[DK] E.N.Domanskii and V.M.Kadets,  On  the  basic  regularizability
of
inverse  operators  (Russian),  Sibirsk.  Mat.  Zh.  29   (1988),
104-108.

[F] V.P.Fonf, Operational bases and  generalized  summation  bases,
Dokl. Akad. Nauk Ukrain. SSR, Ser. A (1986), no. 11,  p.  16-18.
(Russian, Ukrainian).

[G] L.V.Gladun, On Banach spaces the conjugates  of  which  contain
norming nonquasibasic subspaces, Engl. transl.:  Siberian  Math.
J. 28 (1987), 220-223.

[GP] L.V.Gladun and A.N.Plichko, On  norming  and  strongly  norming
subspaces of a dual Banach space, Ukr.  Matem.  Zh.  36  (1984),
427-433. (Russian).

[J] R.C.James,  Separable  conjugate  space,  Pacif.  J.  Math.  10
(1960), 563-571.

[JRZ] W.B.Johnson,  H.P.Rosenthal  and  M.Zippin,  On  bases,  finite
dimensional  decompositions  and  weaker  structures  in  Banach
spaces, Israel J. Math. 9 (1971), 488-506.

[K] M.I.Kadets, Non-linear operatorial bases  in  a  Banach  space,
Teor. Funktsii, Funktsion. Anal. i Prilozhen. 2 (1966), 128-130.
(Russian).

[L] J.Lindenstrauss, On James' paper "Separable conjugate spaces",
Israel J. Math. 9 (1971), 279-284.

[LT1] J.Lindenstrauss and L.Tzafriri, Classical  Banach  spaces,
Lecture
Notes in Math., 338 (1973).

[LT2] J.Lindenstrauss and L.Tzafriri,  Classical  Banach  spaces  I.
Sequence spaces, Berlin, Springer, 1977.

[MR] B.Maurey and H.P.Rosenthal, Normalized weakly null sequence with
no
unconditional subsequence, Stud. Math. 61 (1977), 77-98.

[MP] L.D.Menikhes  and  A.N.Plichko,  Conditions  for  linear  and
finite-dimensio\-nal regularizability of linear inverse  problems,
Dokl. Akad. Nauk SSSR, 241 (1978), 1027-1030. (Russian).

[O1] M.I.Ostrovskii, Total property of bounded approximation,  Sib.
Matem. Zh. 30 (1989), no. 3, 180-181,  Engl.  transl.:  Siberian
Math. J. 30 (1989), no. 3, 488-489.

[O2] M.I.Ostrovskii, Basic and quasibasic  subspaces  in  dual  Banach
spaces, Math. Notes 47 (1990), 584-588.

[O3] M.I.Ostrovskii, Regularizability of inverse linear operators in
Banach
spaces with bases, Siberian Math. J. 33 (1992), 470--476.

[PP] Yu.I.Petunin and A.N.Plichko, The theory of characteristic  of
subspaces and its  applications,  Kiev,  Vyshcha  Shkola,  1980.
(Russian).

[R] H.P.Rosenthal,  On  the  subspaces  of $L^{p}\ (p>2)$  spanned  by
sequences of independent random variables,  Israel  J.  Math.  8
(1970), 273-303.

[S1] I.Singer, On  Banach  spaces  in  which  every $M$-basis  is  a
generalized  summation  basis,  Banach  Center  Publications,  4
(1979), 237-240.

[S2] I.Singer, Bases in Banach spaces. II, Berlin, Springer, 1981.

[V1] F.S.Vakher, The local  problem  of  existence  of  operatorial
bases in Banach spaces, Sib. Matem.  Zh.,  16  (1975),  853-855.
(Russian).

[V2] F.S.Vakher, Bounded approximation property in separable Banach
spaces, Dokl. Akad. Nauk SSSR 255 (1980), 1301-1306. (Russian).

[VP] F.S.Vakher and A.N.Plichko, The bounded approximation property
and linear finite-dimesional regularizability, Ukr.  Matem.  Zh.
33 (1981), 167-171. (Russian).

[VGP]  V.A.Vinokurov,  L.V.Gladun  and   A.N.Plichko,   On   norming
subspaces in a dual Banach space  and  the  regularizability  of
inverse operators, Izvestiya Vuzov. Ser. matem. (1985),  no.  6,
3-10. (Russian).

[ViP] V.A.Vinokurov and  A.N.Plichko,  On  the  regularizability  of
linear inverse problems by  linear  methods.  Dokl.  Akad.  Nauk
SSSR, 229 (1976), 1037-1040. (Russian).

[Z] M.Zippin, A remark on bases and  reflexivity  in  Banach  spaces,
Isr. J. Math. 6 (1968), 74-79.

\enddocument